\documentclass[10pt]{article}

\usepackage[T1]{fontenc}
\usepackage{lmodern}
\usepackage[margin=1in]{geometry}
\usepackage{amsmath,amssymb,amsthm,mathtools}
\usepackage{microtype}
\usepackage[numbers,sort&compress]{natbib}
\usepackage[colorlinks=true,linkcolor=blue,citecolor=red,urlcolor=blue]{hyperref}

\newtheorem{theorem}{Theorem}[section]
\newtheorem{lemma}[theorem]{Lemma}

\newtheorem{corollary}[theorem]{Corollary}
\theoremstyle{definition}

\newtheorem{conjecture}[theorem]{Conjecture}

\newcommand{\FS}{\mathsf{FS}}
\newcommand{\Star}{\operatorname{Star}}

\title{Local connectivity of friends-and-strangers graphs with a star}

\author{Jiaxiang Yang\\[0.6ex]
\normalsize Center for Combinatorics, Nankai University,\\[-0.2ex]
\normalsize Tianjin 300071, P.R. China\\[0.4ex]
\normalsize Email address: corresponding.yjxmath@163.com}

\date{October 31, 2026}

\begin{document}
\maketitle

\begin{abstract}
Given two finite simple graphs $X$ and $Y$ with the same number of vertices, the vertex set of the friends-and-strangers graph $\FS(X,Y)$ consists of all bijections from the vertices of $X$ to the vertices of $Y$. Two bijections are adjacent exactly when one is obtained from the other by swapping the images of two adjacent vertices of $X$, while leaving all other images unchanged, and the two swapped images are adjacent in $Y$. For $n\ge2$, let $\Star_n$ be the star on $n$ vertices.  Krishnan and Li [Electron. J. Combin. 33 (2026), no. 3, Paper P3.41] conjectured that, for $n\ge3$, whenever $\FS(X,\Star_n)$ is connected, the maximum number of pairwise internally vertex-disjoint paths between any two distinct vertices is equal to the smaller of their degrees. We prove this conjecture. More precisely,  if $\FS(X,\Star_n)$ is connected, then $X$ must be $2$-connected, and the conjectured equality holds for every pair of distinct vertices. We further prove the stronger result that, whenever $X$ is $2$-connected, the same equality holds in every connected component of $\FS(X,\Star_n)$.
\end{abstract}

\noindent\textit{2020 Mathematics subject classification:} Primary 05C40; Secondary 05C38.

\noindent\textit{Keywords:} Friends-and-strangers graph; Local connectivity;  Vertex-disjoint paths; Sliding puzzle.

\section{Introduction}

\noindent
For a finite simple graph $G$, let $V(G)$ and $E(G)$ denote its vertex set and edge set, respectively. For a finite set $A$, let $|A|$ denote its cardinality. Let $X$ and $Y$ be finite simple graphs with $|V(X)|=|V(Y)|$.  The \emph{friends-and-strangers graph} $\FS(X,Y)$ is defined as follows. Its vertex set consists of all bijections $\sigma:V(X)\to V(Y)$. Two bijections $\sigma$ and $\sigma'$ are adjacent in $\FS(X,Y)$ if and only if there is an edge $xy\in E(X)$ such that $\{\sigma(x),\sigma(y)\}\in E(Y)$, $\sigma'(x)=\sigma(y)$, $\sigma'(y)=\sigma(x)$, and $\sigma'(z)=\sigma(z)$ for every $z\in V(X)\setminus\{x,y\}$. We call such a bijection a \emph{configuration}, and such an exchange a \emph{friendly swap}. Friends-and-strangers graphs were introduced by Colin Defant and Noah Kravitz~\cite{DefantKravitz2021}.

A basic problem is to determine when $\FS(X,Y)$ is connected. Alon, Defant, and Kravitz~\cite{AlonDefantKravitz2023} studied this problem from random and extremal points of view, while Bangachev~\cite{Bangachev2022} considered asymmetric versions of related extremal questions. Wang and Chen~\cite{WangChenRandom2023} investigated connectedness for random pairs of graphs. Several papers considered particular choices of the position graph or the friendship graph. These include lollipop graphs and related families~\cite{WangChenLollipop2023}, complete bipartite graphs and other special graphs~\cite{WangLuChen2023}, and complete multipartite graphs~\cite{Zhu2025}. Other developments include the study of diameters of connected components by Jeong~\cite{Jeong2024}, connectivity questions for old and new variants of the model by Milojevi\'c~\cite{Milojevic2024}, and connectedness conditions involving graph parameters by Zhao, You, Lin, and Zhang~\cite{ZhaoYouLinZhang2025}.

More recently, Krishnan and Li~\cite{KrishnanLi2026} studied vertex connectivity of friends-and-strangers graphs. For $n\ge2$,  let $\Star_n$ be the star with vertex set
$[n]=\{1,\ldots,n\}$ and center $n$. In $\FS(X,\Star_n)$, every friendly swap moves the label $n$ along an edge of $X$, so the position occupied by $n$ may be viewed as an empty position.

For $n\ge3$ and connected $X$, Krishnan and Li proved that every connected component of $\FS(X,\Star_n)$ has vertex connectivity equal to its minimum degree. For distinct vertices $a$ and $b$ of a graph $G$, let $\kappa_G(a,b)$ denote the maximum number of pairwise internally vertex-disjoint $a$--$b$ paths. They proposed the following strengthening.

\begin{conjecture}\cite[Conjecture 43]{KrishnanLi2026}
Let $n\ge3$, and suppose that $\FS(X,\Star_n)$ is connected. Then for every two distinct configurations $\sigma$ and $\rho$,
\[
\kappa_{\FS(X,\Star_n)}(\sigma,\rho)
=
\min\bigl\{
d_{\FS(X,\Star_n)}(\sigma),
d_{\FS(X,\Star_n)}(\rho)
\bigr\}.
\]
\end{conjecture}

We prove the following stronger componentwise statement. A graph is called \emph{$2$-connected} if it has at least three vertices, is connected, and has no cut vertex.

\begin{theorem}\label{thm:main}
Let $n\ge3$, let $X$ be a $2$-connected graph on $n$ vertices, and let $H$ be any connected component of $\FS(X,\Star_n)$. For every two distinct configurations $\sigma,\rho\in V(H)$,
\[
\kappa_H(\sigma,\rho)=\min\{d_H(\sigma),d_H(\rho)\}.
\]
\end{theorem}

Theorem~\ref{thm:main} does not require $\FS(X,\Star_n)$ itself to be connected. It gives the desired equality separately in every connected component whenever $X$ is $2$-connected. We also show that, for $n\ge3$, connectedness of $\FS(X,\Star_n)$ forces $X$ to be $2$-connected. Therefore the connected case of Theorem~\ref{thm:main} gives Conjecture 43 of Krishnan and Li. The case $n=2$ is handled directly.

In Section~\ref{sec:prelim}, we fix the notation and collect the preliminary results used in the proof. In Section~\ref{sec:fan}, we construct disjoint routes from a configuration to configurations having a prescribed empty position. In Section~\ref{sec:aux}, we introduce an auxiliary bipartite graph and determine its vertex connectivity. The proof of Theorem~\ref{thm:main} and the deduction of Conjecture 43 are given in Section~\ref{sec:proof}.

\section{Preliminaries}\label{sec:prelim}

For a graph $G$ and a vertex $v\in V(G)$, let $N_G(v)$ be the set of neighbors of $v$ and let $d_G(v)=|N_G(v)|$. We write $\delta(G)$ for the minimum degree of $G$ and $\kappa(G)$ for its vertex connectivity, that is, the minimum number of vertices whose deletion disconnects $G$ or leaves a single vertex.

For a configuration $\sigma\in V(\FS(X,\Star_n))$, we call $\sigma^{-1}(n)$ its empty position. If $H$ is a connected component of $\FS(X,\Star_n)$ and $\sigma\in V(H)$, the friendly swaps from $\sigma$ are exactly the moves of the empty position to a neighboring vertex of $X$. Therefore
\[
d_H(\sigma)=d_X\bigl(\sigma^{-1}(n)\bigr).
\]

We use the usual vertex form of Menger's theorem: if two vertices of a finite graph are nonadjacent, then the maximum number of pairwise internally vertex-disjoint paths between them equals the minimum size of a vertex set separating them.

Let $G$ be a bipartite graph with bipartition $(L,R)$. We call $G$ half vertex-transitive if some subgroup of its automorphism group preserves $L$ and $R$ setwise and acts transitively on $L$ and on $R$. For a bipartite digraph, we use the same definition with
automorphisms preserving arc directions. For a strongly connected digraph $D$, let $\kappa(D)$ denote
the minimum number of vertices whose deletion destroys strong
connectivity, and let $\delta(D)$ be the smaller of its minimum
in-degree and minimum out-degree.

We need the following consequence.

\begin{lemma}\cite[Corollary 3.7]{ChenMengTianLiangLiu2018} \label{lem:Chen}
Let $D$ be a strongly connected half vertex-transitive digraph. Then
\[
\kappa(D)=\delta(D).
\]
\end{lemma}

The next lemma gives the undirected form used later.

\begin{lemma}\label{lem:halftrans}
Let $G$ be a connected half vertex-transitive  graph with $\delta(G)\ge2$. Then
\[
\kappa(G)=\delta(G).
\]
\end{lemma}

\begin{proof}
Replace every edge of $G$ by two opposite arcs and denote the resulting symmetric digraph by $D$. Since $G$ is connected, $D$ is strongly connected. The subgroup witnessing the half vertex-transitivity of $G$ also acts transitively on the two parts of $D$, so $D$ is half vertex-transitive. Moreover, its minimum in-degree and minimum out-degree are both $\delta(G)$, and hence $\delta(D)=\delta(G)$.

By Lemma~\ref{lem:Chen}, $\kappa(D)=\delta(G)$. Therefore deleting fewer than $\delta(G)$ vertices leaves $D$ strongly connected. Its underlying undirected graph is the corresponding vertex-deleted subgraph of $G$, so deleting fewer than $\delta(G)$ vertices leaves $G$ connected. Thus $\kappa(G)\ge\delta(G)$.

For the reverse inequality, choose a vertex $v$ with $d_G(v)=\delta(G)$. If $v$ is not adjacent to every other vertex, deleting $N_G(v)$ isolates $v$ while at least one other vertex remains. If $v$ is adjacent to every other vertex, then $\delta(G)=|V(G)|-1$, and deleting $V(G)\setminus\{v\}$ leaves a single vertex. In either case, $\kappa(G)\le\delta(G)$. Therefore $\kappa(G)=\delta(G)$.
\end{proof}

The next lemma explains why Theorem~\ref{thm:main} contains the connected case of Conjecture 43.

\begin{lemma}\label{lem:2connected}
Let $n\ge3$. If $\FS(X,\Star_n)$ is connected, then $X$ is $2$-connected.
\end{lemma}

\begin{proof}
First, $X$ must be connected. Indeed, the empty position moves only along edges of $X$, so it cannot move between two different connected components of $X$.

Suppose for a contradiction that $X$ has a cut vertex $u$. Let
$
C_1,C_2,\ldots,C_t
$
be the connected components of $X-u$, where $t\ge2$. Choose a configuration $\sigma$ whose empty position is $u$, so $\sigma(u)=n$.

Consider any sequence of friendly swaps that starts at $\sigma$ and ends at a configuration whose empty position is again $u$. Each time the empty position leaves $u$, it enters one component $C_i$. Before it can enter a different component of $X-u$, it must first return to $u$.

Consider one excursion that starts with the empty position at $u$, enters $C_i$, and ends when the empty position first returns to $u$. At the first move, one label from $C_i$ moves to $u$. Until the final move of the excursion, the empty position stays inside $C_i$, so the label left at $u$ cannot move. At the last move, that same label moves from $u$ back into $C_i$. Hence the set of labels occupying $C_i$ is the same immediately before and immediately after the excursion. In particular, no excursion transfers a label from one component of $X-u$ to another.

It follows that whenever a reachable configuration has empty position $u$, the set of labels occupying each $C_i$ is the same as in $\sigma$. Now interchange one noncentral label in $C_1$ with one noncentral label in $C_2$, while keeping the label $n$ at $u$. This gives another configuration of $\FS(X,\Star_n)$, but it cannot be reached from $\sigma$. This contradicts the assumed connectedness of $\FS(X,\Star_n)$.

Therefore $X$ has no cut vertex. Since $n\ge3$ and $X$ is connected, $X$ is $2$-connected.
\end{proof}

\section{Disjoint routes to a fixed empty position}\label{sec:fan}

Let $X$ be $2$-connected and $H$ a fixed connected component of $\FS(X,\Star_n)$. Because $X$ is connected, the empty position can be moved along a path of $X$ from any vertex to any other vertex. Therefore $H$ contains a configuration whose empty position is any prescribed vertex of $X$.

Fix $r\in V(X)$ and define
\[
F_r=\{\tau\in V(H):\tau^{-1}(n)=r\}.
\]
Thus $F_r$ is the set of configurations in $H$ whose empty position is $r$. We shall call $F_r$ the \emph{$r$-fiber}.

\begin{lemma}\label{lem:fan}
Let $\sigma\in V(H)\setminus F_r$. There are $d_H(\sigma)$ paths from $\sigma$ to $F_r$ such that any two of them meet only at their common initial vertex $\sigma$. Their endpoints in $F_r$ are distinct.
\end{lemma}

\begin{proof}
Let $u=\sigma^{-1}(n)$ be the empty position of $\sigma$. Since $\sigma\notin F_r$, we have $u\ne r$. Because $X$ is $2$-connected, deleting $u$ leaves a connected graph $X-u$.

Take a neighbor $v\in N_X(u)$. If $v=r$, use the edge $ur$. Otherwise, since $X-u$ is connected, choose a simple $v$--$r$ path in $X-u$ and prepend the edge $uv$. In either case, we obtain a simple $u$--$r$ path in $X$ whose first edge is $uv$ and which does not return to $u$ after the first step.

Starting from $\sigma$, move the empty position along this $u$--$r$ path. Every move is a friendly swap because the empty label is the center of $\Star_n$. Since the empty position visits distinct vertices along the chosen path, the resulting configurations are distinct. Thus they form a path in $H$ from $\sigma$ to a configuration in $F_r$, and this path meets $F_r$ only at its endpoint. We construct one such path for every $v\in N_X(u)$.

It remains to prove that these configuration paths are disjoint after leaving $\sigma$. Take two distinct neighbors $v,w\in N_X(u)$. On the route whose first move is across $uv$, the label $\sigma(v)$ is moved to $u$ at the first step. The empty position never returns to $u$, so the label at $u$ remains $\sigma(v)$ at every later configuration on that route. Similarly, every configuration after the first step on the route beginning with $uw$ has the label $\sigma(w)$ at $u$.

The bijection $\sigma$ assigns different labels to $v$ and $w$, so $\sigma(v)\ne\sigma(w)$. Therefore no configuration after the first step of the first route can equal a configuration after the first step of the second route. The same argument shows that their final configurations in $F_r$ are distinct.

Thus the constructed routes are pairwise disjoint outside $\sigma$. Their number is
$
|N_X(u)|=d_X(u)=d_H(\sigma),
$
where the last equality follows from the degree identity.
\end{proof}

\section{The auxiliary bipartite graph}\label{sec:aux}

Let $H-F_r$ be the graph obtained from $H$ by deleting every vertex in $F_r$ and all incident edges. Let $\mathcal C_r$ be the set of connected components of $H-F_r$.

Define a bipartite graph $B_r$ with bipartition $(F_r,\mathcal C_r)$. For $\tau\in F_r$ and $C\in\mathcal C_r$, join $\tau$ to $C$ in $B_r$ if $\tau$ has at least one neighbor in $C$ in the graph $H$.

\begin{lemma}\label{lem:aux}
For every $r\in V(X)$,
\[
\kappa(B_r)=d_X(r).
\]
\end{lemma}

\begin{proof}
We verify connectedness, minimum degree, and half vertex-transitivity separately.

First, $B_r$ is connected.
Contract every connected component of $H-F_r$ to one vertex, delete the resulting loops, and then suppress multiple edges. No two configurations in $F_r$ are adjacent in $H$, because a friendly swap necessarily moves the empty position away from $r$. Hence $F_r$ is an independent set. The resulting simple bipartite graph is exactly $B_r$. Since $H$ is connected, the contraction is connected, so $B_r$ is connected.

Second, vertices in $F_r$ have degree $d_X(r)$.
Fix $\tau\in F_r$. For every $v\in N_X(r)$, moving the empty position across the edge $rv$ gives one neighbor of $\tau$ in $H$. These are all the neighbors of $\tau$, so there are exactly $d_X(r)$ of them.

We claim that the neighbors obtained from two different vertices $v,w\in N_X(r)$ lie in different components of $H-F_r$. After the move from $r$ to $v$, the label $\tau(v)$ is left at $r$. Along any path in $H-F_r$, the empty position never visits $r$, so the label at $r$ cannot change. Similarly, after the move from $r$ to $w$, the label at $r$ is $\tau(w)$. Since $\tau(v)\ne\tau(w)$, the two resulting configurations cannot lie in the same component of $H-F_r$.

Thus the $d_X(r)$ neighbors of $\tau$ in $H$ meet $d_X(r)$ different component vertices of $B_r$. Therefore
$
d_{B_r}(\tau)=d_X(r).
$

Third, vertices in $\mathcal C_r$ have degree at least $d_X(r)$.
Fix $C\in\mathcal C_r$ and choose a configuration in $C$. Its empty position is different from $r$. Since $X$ is $2$-connected, the graph $X-r$ is connected. Therefore, while avoiding $r$, the empty position can be moved to any prescribed vertex of $X-r$. Every configuration on such a route remains in $H-F_r$, and hence remains inside the same component $C$.

For each $v\in N_X(r)$, choose $\tau_v\in C$ whose empty position is $v$. Move the empty position once across $vr$ and call the resulting configuration $\widehat\tau_v\in F_r$. Then $C$ is adjacent in $B_r$ to every $\widehat\tau_v$.

These configurations are distinct. Suppose that $\widehat\tau_v=\widehat\tau_w$ for two distinct neighbors $v,w$ of $r$. Then this common configuration in $F_r$ would have two neighbors, namely $\tau_v$ and $\tau_w$, lying in the same component $C$ of $H-F_r$. This contradicts the preceding degree argument for vertices in $F_r$, which showed that different first moves from a configuration in $F_r$ lead to different components of $H-F_r$. Hence
$
d_{B_r}(C)\ge d_X(r).
$
Combining this with the preceding equality gives
$
\delta(B_r)=d_X(r).
$

Finally, $B_r$ is half vertex-transitive.
Let $\Gamma_H$ be the stabilizer of $H$ in the group of permutations of the leaf labels. Explicitly, $\Gamma_H$ consists of the permutations $\varphi$ of $[n]$ such that $\varphi(n)=n$ and the induced map $\eta\mapsto\varphi\circ\eta$ sends $H$ to itself. It is a subgroup under composition. Since every permutation fixing $n$ is an automorphism of $\Star_n$, each element of $\Gamma_H$ induces an automorphism of $H$ and preserves the empty position.

Take $\tau,\pi\in F_r$. Define the label permutation $\varphi=\pi\circ\tau^{-1}$. Because $\tau(r)=\pi(r)=n$, we have $\varphi(n)=n$, so $\varphi$ is an automorphism of $\Star_n$. Moreover,
$
\varphi\circ\tau=\pi.
$
The map $\eta\mapsto\varphi\circ\eta$ is an automorphism of $\FS(X,\Star_n)$ and therefore sends connected components to connected components. Since both $\tau$ and $\pi$ lie in $H$, it sends $H$ to itself. Hence $\varphi\in\Gamma_H$.  This proves that $\Gamma_H$ acts transitively on $F_r$.

Now take $C,D\in\mathcal C_r$ and fix one vertex $v\in V(X)\setminus\{r\}$. As observed above, since $X-r$ is connected, both $C$ and $D$ contain a configuration whose empty position is $v$. Choose $\tau\in C$ and $\pi\in D$ with empty position $v$, and again set $\varphi=\pi\circ\tau^{-1}$. Since $\tau(v)=\pi(v)=n$, we have
$
\varphi(n)=\pi\bigl(\tau^{-1}(n)\bigr)=\pi(v)=n.
$
Thus $\varphi$ fixes the center label of $\Star_n$. Moreover, $\varphi\circ\tau=\pi$, and both $\tau$ and $\pi$ lie in $H$. Therefore the automorphism induced by $\varphi$ sends the connected component $H$ to itself, so $\varphi\in\Gamma_H$. Since $\varphi$ preserves the empty position, it preserves $F_r$ and therefore maps connected components of $H-F_r$ to connected components of $H-F_r$. Because it sends a vertex of $C$ to a vertex of $D$, it sends $C$ to $D$. Thus $\Gamma_H$ is also transitive on $\mathcal C_r$.

Every element of $\Gamma_H$ preserves adjacency in $H$ and preserves $F_r$. It therefore permutes the components of $H-F_r$ and preserves the incidence relation defining $B_r$. Hence $\Gamma_H$ induces a subgroup of $\operatorname{Aut}(B_r)$. The preceding two paragraphs show that this subgroup is transitive on both sides of the bipartition $(F_r,\mathcal C_r)$. Therefore $B_r$ is half vertex-transitive.

Since $X$ is $2$-connected, every vertex of $X$ has degree at least $2$. In particular,
$
d_X(r)\ge2.
$
Together with $\delta(B_r)=d_X(r)$, this verifies the minimum-degree hypothesis of Lemma~\ref{lem:halftrans}. Therefore
$
\kappa(B_r)=\delta(B_r)=d_X(r).
$
\end{proof}

We next transfer the connectivity of $B_r$ back to the original component $H$. The statement allows both deleted vertices and deleted edges. The edge deletion is needed later when the two endpoints under consideration are adjacent.

For a vertex set $S\subseteq V(H)$ and an edge set $E_0\subseteq E(H)$, write $H-S-E_0$ for the graph obtained by deleting all vertices in $S$ together with their incident edges, and then deleting all edges in $E_0$ whose endpoints remain.

\begin{lemma}\label{lem:fiberfault}
Let $S\subseteq V(H)$ and $E_0\subseteq E(H)$. If
$
|S|+|E_0|<d_X(r),
$
then every two configurations in $F_r\setminus S$ are connected in $H-S-E_0$.
\end{lemma}

\begin{proof}
We delete vertices from the auxiliary graph $B_r$ so that any path remaining in $B_r$ can be lifted to a path of $H$ avoiding every element of $S\cup E_0$.

For each configuration in $S\cap F_r$, delete the corresponding vertex on the $F_r$ side of $B_r$. For each configuration $x\in S\setminus F_r$, delete the vertex of $B_r$ corresponding to the unique component of $H-F_r$ that contains $x$.

Now consider an edge $e\in E_0$. The set $F_r$ is independent, so at least one endpoint of $e$ lies outside $F_r$. If exactly one endpoint lies outside $F_r$, delete the component vertex of $B_r$ containing that endpoint. If both endpoints lie outside $F_r$, then the two endpoints are adjacent in $H-F_r$ and therefore lie in the same component of $H-F_r$; delete that component vertex. Thus each edge of $E_0$ causes the deletion of at most one vertex from the $\mathcal C_r$ side of $B_r$.

Altogether, at most $|S|+|E_0|$ vertices are deleted from $B_r$. By Lemma~\ref{lem:aux}, this number is smaller than $\kappa(B_r)$. Hence the remaining auxiliary graph is connected.

Take $\alpha,\beta\in F_r\setminus S$ and choose an $\alpha$--$\beta$ path in the remaining auxiliary graph. Consider a length-two segment
$
\alpha'-C-\beta'
$
of this auxiliary path, where $C\in\mathcal C_r$. By the definition of $B_r$, there is an edge $\alpha'a$ of $H$ with $a\in C$, and there is an edge $b\beta'$ of $H$ with $b\in C$.

The component vertex $C$ was not deleted from $B_r$. Therefore no vertex of $C$ belongs to $S$. Moreover, no edge of $E_0$ has an endpoint in $C$, because otherwise the construction above would have deleted the vertex $C$ from $B_r$. Hence every edge inside $C$ remains present. The vertices $\alpha'$ and $\beta'$ lie on the surviving auxiliary path, so neither belongs to $S$; and the chosen boundary edges $\alpha'a$ and $b\beta'$ are not in $E_0$, since each has an endpoint in $C$. Thus both boundary edges survive in $H-S-E_0$. Since $C$ is connected and all of its vertices and internal edges survive, there is an $a$--$b$ path inside $C$ that also survives.

Replace every segment $\alpha'-C-\beta'$ of the auxiliary path by a simple path
$
\alpha' a\;\cdots\; b\beta'
$
through $C$. Since the auxiliary path visits each component vertex at most once and distinct components of $H-F_r$ are vertex-disjoint, the interiors of these replacement paths are pairwise disjoint. Hence these replacements yield an $\alpha$--$\beta$ path in $H-S-E_0$. Therefore all vertices of $F_r\setminus S$ lie in one connected component of $H-S-E_0$.
\end{proof}

\section{Proof of the main theorem}\label{sec:proof}

\begin{proof}[Proof of Theorem~\ref{thm:main}]
Take two distinct configurations $\sigma,\rho\in V(H)$. Interchanging their names if necessary, assume
$
d_H(\sigma)\le d_H(\rho).
$
Let $r=\rho^{-1}(n)$ be the empty position of $\rho$. Then $\rho\in F_r$. By the degree identity from Section~\ref{sec:prelim},
$
d_X(r)=d_H(\rho)\ge d_H(\sigma).
$

We first establish the fault-tolerance statement that drives the proof. Let
$
S\subseteq V(H)\setminus\{\sigma,\rho\}
$
and let $E_0\subseteq E(H)$ satisfy
$
|S|+|E_0|<d_H(\sigma).
$
We claim that $\sigma$ and $\rho$ remain connected in $H-S-E_0$.

Suppose first that $\sigma\in F_r$. Then both $\sigma$ and $\rho$ belong to $F_r\setminus S$. Since
$
|S|+|E_0|<d_H(\sigma)\le d_X(r),
$
Lemma~\ref{lem:fiberfault} gives a $\sigma$--$\rho$ path in $H-S-E_0$.

Now suppose that $\sigma\notin F_r$. By Lemma~\ref{lem:fan}, there are exactly $d_H(\sigma)$ routes from $\sigma$ to $F_r$ that are pairwise disjoint outside $\sigma$ and have distinct endpoints in $F_r$.

A deleted vertex from $S$ lies on at most one of these routes because the routes are vertex-disjoint outside $\sigma$, and $\sigma\notin S$. A deleted edge from $E_0$ also lies on at most one route. Indeed, the routes have distinct first edges out of $\sigma$, and after the first step they are vertex-disjoint. Hence the faults in $S\cup E_0$ destroy at most $|S|+|E_0|$ of the $d_H(\sigma)$ routes.

Since $|S|+|E_0|<d_H(\sigma)$, at least one route survives. Let $\alpha$ be its endpoint in $F_r$. The endpoint $\alpha$ is not in $S$, because otherwise the surviving route would contain a deleted vertex. Thus $\alpha,\rho\in F_r\setminus S$. Again,
\[
|S|+|E_0|<d_H(\sigma)\le d_X(r),
\]
so Lemma~\ref{lem:fiberfault} connects $\alpha$ to $\rho$ inside $H-S-E_0$. Concatenating the surviving $\sigma$--$\alpha$ route with this $\alpha$--$\rho$ path gives a $\sigma$--$\rho$ walk in $H-S-E_0$. Removing closed repeated sections from this walk yields a $\sigma$--$\rho$ path. This proves the fault-tolerance claim.

We now convert this claim into internally vertex-disjoint paths.

Assume first that $\sigma$ and $\rho$ are not adjacent in $H$. Take $E_0=\varnothing$. The fault-tolerance claim shows that deleting any set of at most $d_H(\sigma)-1$ vertices other than $\sigma$ and $\rho$ does not separate them. By Menger's theorem, there are at least $d_H(\sigma)$ pairwise internally vertex-disjoint $\sigma$--$\rho$ paths.

Now assume that $\sigma$ and $\rho$ are adjacent. Since $X$ is $2$-connected, every vertex of $X$ has degree at least $2$. Thus $d_H(\sigma)=d_X(\sigma^{-1}(n))\ge2$. Delete the edge $\sigma\rho$. If
$
S\subseteq V(H)\setminus\{\sigma,\rho\}
$
has at most $d_H(\sigma)-2$ vertices, then
$
|S|+1<d_H(\sigma).
$
Applying the fault-tolerance claim with $E_0=\{\sigma\rho\}$ shows that $\sigma$ and $\rho$ remain connected in $H-\sigma\rho-S$. Thus no set of at most $d_H(\sigma)-2$ vertices separates $\sigma$ from $\rho$ in the graph $H-\sigma\rho$. They are nonadjacent in this graph, so Menger's theorem gives at least $d_H(\sigma)-1$ pairwise internally vertex-disjoint $\sigma$--$\rho$ paths in $H-\sigma\rho$. Adding the one-edge path $\sigma\rho$ gives $d_H(\sigma)$ pairwise internally vertex-disjoint paths in $H$.
 We have proved
$
\kappa_H(\sigma,\rho)\ge d_H(\sigma).
$
For the reverse inequality, consider any family of pairwise internally vertex-disjoint $\sigma$--$\rho$ paths. Apart from the possible one-edge path $\sigma\rho$, different paths must leave $\sigma$ through different neighbors and must enter $\rho$ through different neighbors. Therefore the number of paths is at most both endpoint degrees. Hence
\[
\kappa_H(\sigma,\rho)
\le \min\{d_H(\sigma),d_H(\rho)\}=d_H(\sigma).
\]
The upper and lower bounds are equal, which proves the theorem.
\end{proof}

\begin{corollary}\label{cor:KL}
Let $n\ge2$, and suppose that $\FS(X,\Star_n)$ is connected. Then for every two distinct configurations $\sigma$ and $\rho$,
\[
\kappa_{\FS(X,\Star_n)}(\sigma,\rho)
=
\min\bigl\{d_{\FS(X,\Star_n)}(\sigma),d_{\FS(X,\Star_n)}(\rho)\bigr\}.
\]
\end{corollary}

\begin{proof}
If $n=2$, connectedness forces $X=K_2$. In that case $\Star_2=K_2$ and $\FS(K_2,K_2)\cong K_2$, so the unique pair of distinct configurations is joined by exactly one path and both endpoint degrees are $1$.

Now assume $n\ge3$. By Lemma~\ref{lem:2connected}, the connectedness of $\FS(X,\Star_n)$ implies that $X$ is $2$-connected. The graph $\FS(X,\Star_n)$ is then the unique connected component to which Theorem~\ref{thm:main} applies.
\end{proof}

\section*{Acknowledgments}
This work was supported by the Fundamental and Interdisciplinary Disciplines Breakthrough Plan of the Ministry of Education of China (JYB2025XDXM207) and  the National Natural Science Foundation of China (No. 12631013).

\section*{Declaration of AI assistance}
During the preparation of this manuscript, the author used ChatGPT to assist in exploring possible approaches, checking technical details, and improving the exposition. All mathematical arguments, computations, proofs, and results were independently verified by the author. The author takes full responsibility for the content of this manuscript.

\end{document}